\documentclass[12]{article}

\usepackage[T2A]{fontenc}
\usepackage[cp1251]{inputenc}
\usepackage[russian,english]{babel}
\usepackage[tbtags]{amsmath}
\usepackage{amssymb,amsfonts}
\usepackage{mathrsfs}
\usepackage[12pt]{extsizes}
\begin{document}
\date{}

\author{D. O. Shatskov}

\title{On the irrationality  measure function in average.}

\maketitle

\begin{abstract}
We consider the value $I_\alpha(t)=\int\limits_1^t\psi_\alpha
(\xi)d\xi$, where $\psi_\alpha(t)=\min\limits_{1\leqslant
q\leqslant t}||q\alpha||$. We prove that for almost all $\alpha$
one has $\lim\limits_{t\rightarrow+\infty}\frac{I_\alpha(t)}{\ln
t} = \frac{6\log 2}{\pi^2}$. It is proved that there  exist
algebraically independent numbers $\alpha$ and $\beta $ such that
the difference $ I_\alpha (t) - I_\beta (t)$ tends to infinity
when $t \to +\infty.$

Bibliography: 10 titles.
\end{abstract}

\footnotetext{ Research is supported by the grant RFBR No.
12-01-00681-a.} \vskip+1.0cm

Let  $\alpha$ be  a real irrational number. We consider the
funcion
$$
\psi_\alpha(t)=\min\limits_{1\leqslant q\leqslant t}||q\alpha||
$$
(here $q$ is an integer number and $\|\cdot\|$ stands for the
distance to the nearest integer). For $0<\alpha<1$  we consider
the continued fraction expansion
$$
\alpha = [a_0;a_1,a_2,\ldots] = a_0
+\frac{1}{a_1+\displaystyle{\frac{1}{a_2+\cdots}}}.
$$
Let us define the convergents
$$
\frac{p_\nu}{q_{\nu}} = [a_0;,a_1,a_2,\ldots,a_\nu].
$$
For $t$ from the range $q_\nu\leqslant t<q_{\nu+1}$ one has
$$
\psi_\alpha (t) = ||q_{\nu} \alpha ||= |q_\nu \alpha - p_\nu |.
\eqno (1)
$$

Let us recall the definition of  {\it Lagrange specrtum } (see
\cite{Cusick}):
$$
\mathbb{L} = \{ \lambda \in \mathbb{R}:\,\,\, \exists \,\alpha
\,\,\, \liminf_{t\to +\infty} t\psi_\alpha (t) = \lambda\},
$$
and {\it Dirichlet  spectrum } (see \cite{Ivanov}):
$$
\mathbb{D} = \{ \lambda \in \mathbb{R}:\,\,\, \exists \,\alpha
\,\,\, \limsup_{t\to +\infty} t\psi_\alpha (t) = \lambda\}.
$$
There is  a lot of results related to the structure of the sets
$\mathbb{L},\,\mathbb{D}$. In the present paper we study the
integral
$$
I_\alpha(t)=\int\limits_1^t\psi_\alpha (\xi)d\xi.
$$
As $0<t\psi_\alpha(t)<1$ for  $t\geqslant1$, we have
$I_\alpha(t)<\ln t$.

Given $t>1$ we define  $N=N(\alpha,t)$ by the conditions
$$
q_N^{}\leqslant t <q_{N+1}^{}.\eqno(2)
$$

It is clear that $N\rightarrow\infty$ when $t\rightarrow\infty$.

The main results of the present paper are as follows.

{\bf Theorem 1.} {\it For almost all $($in sense of Lebesgue
measure$)$ numbers $\alpha$ one has}

$ 1)\,
\lim\limits_{t\rightarrow\infty}\frac{I_\alpha(t)}{N(\alpha,t)}=\frac{1}{2},$

$ 2)\, \lim\limits_{t\rightarrow\infty}\frac{I_\alpha(t)}{\ln
t}=\frac{6\ln2}{\pi^2}. $\\

In the next two  theorem we have calculated the extremal values
for the quantity $\frac{I_\alpha(t)}{N(\alpha,t)}$.

{\bf Theorem 2.} {\it For every irrational number $\alpha\in(0;1)$
one has}

$1)\,
\limsup\limits_{t\rightarrow\infty}\frac{I_\alpha(t)}{N(\alpha,t)}\leqslant1,$

$2)\,
\liminf\limits_{t\rightarrow\infty}\frac{I_\alpha(t)}{N(\alpha,t)}\geqslant\frac{1}{2}-\frac{\sqrt{5}}{10}.$

The bounds from Theorem 2 are optimal ones.

{\bf Theorem 3.} {\it Given}
$d\in\left[\frac{1}{2}-\frac{\sqrt{5}}{10};1\right]$ {\it there
exists $\alpha$, such that}
$$
\lim\limits_{t\rightarrow\infty}\frac{I_\alpha(t)}{N(\alpha,t)}=d.
$$
In the paper \cite{MosKan} it is  proved that  for any two real
numbers $\alpha, \beta$, such that $\alpha \pm \beta \not\in
\mathbb{Z}$, the difference function
$$
\psi_\alpha (t) - \psi_\beta (t) \eqno(3)
$$
changes its sign
infinitely many often  as $t\rightarrow\infty$.

It is easy  to show that for $\tau = \frac{\sqrt{5}-1}{2}$ one has
$$
\lim\limits_{t\rightarrow\infty}\frac{I_\tau(t)}{\ln t}=
\left(\frac{1}{2}-\frac{\sqrt{5}}{10}\right):\ln\left(\frac{\sqrt{5}+1}{2}\right).
$$
So one can easily find two algebraically independent
 real  numbers $\alpha$ and $\beta $ such that
the limits
$$
\lim\limits_{t\rightarrow\infty}\frac{I_\alpha(t)}{\ln t} \,\,\,\,
\text{and}
\,\,\,\,\lim\limits_{t\rightarrow\infty}\frac{I_\beta(t)}{\ln t}
$$
are different. We see that  it may happen that the difference $
I_\alpha (t) - I_\beta (t)$ tends to infinity as $t \to +\infty.$
This shows that there is no general oscillating property  for the
difference function
 $I_\alpha (t) - I_\beta (t)$.

For the integrals under consideration we have the following
result.

{\bf Theorem 4.} {\it Let $q_0,q_1,\ldots,q_n,\ldots$ and
$r_0,r_1,\ldots,r_n,\ldots$ be the denominators of the convergents
fractions for $\alpha$ and $\beta$, respectively. Then for almost
all $($in sense of Lebesgue measure) pairs $(\alpha,\beta)$ in
$[0,1]^2$ one has}
$$
\liminf\limits_{n\rightarrow\infty}\left|\int\limits_1^{q_n}\psi_\alpha(\xi)d\xi-
\int\limits_1^{r_n}\psi_\beta(\xi)d\xi\right|<+\infty.
$$
\selectlanguage{russian}
\newpage
{\bf О среднем значении меры иррациональности вещественных чисел.}
\begin{center}
{\bf Д. О. Шацков}
\end{center}

\footnotetext{Работа выполнена при поддержке РФФИ №
12-01-00681-а.}

{\bf 1. Введение и формулировки результатов.}

Для действительного $\alpha$ рассмотрим функцию
$$
\psi_\alpha(t)=\min\limits_{1\leqslant q\leqslant t}||q\alpha||
$$
(здесь минимум берется по целым $q$ и $\|\cdot\|$ обозначает
расстояние до ближайшего целого). Если рассмотреть разложение
числа $\alpha$ в обыкновенную цепную дробь
$$
\alpha = [a_0;a_1,a_2,\ldots] = a_0
+\frac{1}{a_1+\displaystyle{\frac{1}{a_2+\cdots}}}
$$
и через
$$
\frac{p_\nu}{q_{\nu}} = [a_0;,a_1,a_2,\ldots,a_\nu]
$$
обозначать подходящие к $\alpha$ дроби, то при $q_\nu\leqslant
t<q_{\nu+1}$ будет иметь место
\begin{equation}\label{psi}
\psi_\alpha (t) = ||q_{\nu} \alpha ||= |q_\nu \alpha - p_\nu |.
\end{equation}

Известны многочисленные исследования (см. \cite{Cusick}) {\it
спектра Лагранжа}
$$
\mathbb{L} = \{ \lambda \in \mathbb{R}:\,\,\, \exists \,\alpha
\,\,\, \liminf_{t\to +\infty} t\psi_\alpha (t) = \lambda\}
$$
и (см. \cite{Ivanov}) {\it спектра Дирихле}
$$
\mathbb{D} = \{ \lambda \in \mathbb{R}:\,\,\, \exists \,\alpha
\,\,\, \limsup_{t\to +\infty} t\psi_\alpha (t) = \lambda\}.
$$

Предметом нашего изучения будет интеграл
$$
I_\alpha(t)=\int\limits_1^t\psi_\alpha (\xi)d\xi.
$$
Поскольку $0<t\psi_\alpha(t)<1$ для любого $t\geqslant1$, то сразу
видим, что $I_\alpha(t)<\ln t$.

Через $N=N(\alpha,t)$ мы обозначим величину, задаваемую условием
\begin{equation}
 q_N^{}\leqslant t <q_{N+1}^{}.\label{aa}
\end{equation}

Ясно, что $N\rightarrow\infty$ при $t\rightarrow\infty$.

В настоящей работе мы докажем следующий метрический результат.

{\bf Теорема 1.} {\it Для почти всех $($в смысле меры Лебега$)$
чисел $\alpha$ выполняются равенства}

$ 1)\,
\lim\limits_{t\rightarrow\infty}\frac{I_\alpha(t)}{N(\alpha,t)}=\frac{1}{2},$

$ 2)\, \lim\limits_{t\rightarrow\infty}\frac{I_\alpha(t)}{\ln
t}=\frac{6\ln2}{\pi^2}. $\\

Помимо метрического результата, мы докажем утверждение об
экстремальных значениях величины
$\frac{I_\alpha(t)}{N(\alpha,t)}$.

{\bf Теорема 2.} {\it Для любого иррационального $\alpha\in(0;1)$
выполнены неравенства}

$1)\,
\limsup\limits_{t\rightarrow\infty}\frac{I_\alpha(t)}{N(\alpha,t)}\leqslant1,$

$2)\,
\liminf\limits_{t\rightarrow\infty}\frac{I_\alpha(t)}{N(\alpha,t)}\geqslant\frac{1}{2}-\frac{\sqrt{5}}{10}.$

Оценки, приводимые в теореме 2 точны. Более того, имеет место

{\bf Теорема 3.} {\it Для любого}
$d\in\left[\frac{1}{2}-\frac{\sqrt{5}}{10};1\right]$ {\it существует
$\alpha$, такое что}
$$
\lim\limits_{t\rightarrow\infty}\frac{I_\alpha(t)}{N(\alpha,t)}=d.
$$

В недавней работе \cite{MosKan} доказано, что для двух
вещественных чисел $\alpha, \beta$, таких что  $\alpha \pm \beta
\not\in \mathbb{Z}$ разность
\begin{equation}\label{1}
\psi_\alpha (t) - \psi_\beta (t)
\end{equation}
бесконечно много раз  меняет знак при $t\rightarrow\infty$.

Из формул  (\ref{odnzam}), (\ref{ogrS}), (\ref{DogrI}) и леммы 2
настоящей статьи можно получить, что для $\tau =
\frac{\sqrt{5}-1}{2}$ выполнено
$$
\lim\limits_{t\rightarrow\infty}\frac{I_\tau(t)}{\ln t}=
\left(\frac{1}{2}-\frac{\sqrt{5}}{10}\right):\ln\left(\frac{\sqrt{5}+1}{2}\right).
$$
Следовательно существуют алгебраически независимые  $\alpha$ и
$\beta $, такие что пределы
$$
\lim\limits_{t\rightarrow\infty}\frac{I_\alpha(t)}{\ln t} \,\,\,\,
\text{и}
\,\,\,\,\lim\limits_{t\rightarrow\infty}\frac{I_\beta(t)}{\ln t}
$$
различны, и, следовательно разность $ I_\alpha (t) - I_\beta (t)$ может стремится к бесконечности при $t \to +\infty.$
Таким образом, аналог теоремы об
осцилляции разности (\ref{1}) для интегралов $I_\alpha (t),
I_\beta (t)$, вообще-говоря, не имеет места. В настоящей статье мы
доказываем несколько более слабый результат для рассматриваемых
нами интегралов.

{\bf Теорема 4.} {\it Для почти всех $($в смысле меры Лебега на
$[0,1]^2)$ пар $(\alpha,\beta)$ верно следующее утверждение. Пусть
$q_0,q_1,\ldots,q_n,\ldots$ и $r_0,r_1,\ldots,r_n,\ldots$ суть
знаменатели подходящих дробей для $\alpha$ и $\beta$
соответственно. Тогда
$$
\liminf\limits_{n\rightarrow\infty}\left|\int\limits_1^{q_n}\psi_\alpha(\xi)d\xi-
\int\limits_1^{r_n}\psi_\beta(\xi)d\xi\right|<+\infty.
$$
}

Сначала в пунктах 2 и 3 мы приводим доказательство теоремы 2.
Затем в пункте 4 мы доказываем теорему 3. В пункте 5 мы приводим
необходимые сведения из эргодической теории, а затем в пунктах 6 и
7 доказываем теоремы 1 и 4.

{\bf 2. Формулы с подходящими дробями.}

Нам понадобится формула для погрешности приближения числа $\alpha$
его подходящей дробью $\frac{p_n}{q_n}$, которая имеет вид
\begin{equation}\label{formula4}
||q_\nu \alpha || = \frac{1}{q_\nu\alpha_{\nu+1} +q_{\nu-1}},
\end{equation}
или
\begin{equation}\label{formula}
q_\nu ||q_\nu \alpha || = \frac{1}{\alpha_{\nu+1} +\alpha^*_\nu},
\end{equation}
где
$$
\alpha_{\nu+1} = [a_{\nu+1};a_{\nu+2},...],\,\,\,\, \alpha_\nu^* =
[0;a_\nu,a_{\nu-1},...,a_1]=\frac{q_{\nu-1}}{q_\nu}.
$$

Формулы (\ref{psi}) и (\ref{formula4}) позволяют записать интеграл
$I_\alpha (t)$ в виде
\begin{equation}\label{vid}
I_\alpha(t) =
\sum\limits_{\nu=1}^N\frac{q_\nu-q_{\nu-1}}{q_{\nu-1}\alpha_\nu+q_{\nu-2}}+
\frac{t-q_N^{}}{q_N^{}\alpha_{N+1}^{}+q_{N-1}^{}}, \end{equation}
где $N$ определено в (\ref{aa}). Эту же формулу можно записать
по-другому. Введем новое обозначение для слагаемых  из формулы
(\ref{vid}) и преобразуем их:
\begin{equation}\label{ds}
S_\nu(\alpha) =\frac{q_\nu-q_{\nu-1}}
{q_{\nu-1}\left(a_\nu+\frac{1}{\alpha_{\nu+1}}\right)+q_{\nu-2}}=
\frac{q_\nu-q_{\nu-1}}{q_\nu+\frac{q_{\nu-1}}{\alpha_{\nu+1}}}=
\frac{(q_\nu-q_{\nu-1})\alpha_{\nu+1}}{q_\nu\alpha_{\nu+1}+q_{\nu-1}}=
\frac{(1-\alpha_\nu^\ast)\alpha_{\nu+1}}{\alpha_{\nu+1}+\alpha_\nu^\ast}.
\end{equation}
Для последнего слагаемого  в формуле (\ref{vid}) аналогично
получаем
$$
A_{N+1}(\alpha, t)=
\frac{t-q_N}{q_N\alpha_{N+1}+q_{N-1}}=\frac{\left(\frac{t}{q_{N+1}}-\alpha_{N+1}^\ast\right)\alpha_{N+2}}
{\alpha_{N+2}+\alpha_{N+1}^\ast}.
$$
Теперь формулу (\ref{vid}) можно записать
\begin{equation}\label{s}
I_\alpha (t) = G_N(\alpha) + A_{N+1} (\alpha, t), \,\,\,\,
\end{equation}
где
\begin{equation}\label{summa}
G_N(\alpha) = \sum_{\nu = 1}^N S_\nu (\alpha).
\end{equation}
Поскольку $\alpha_\nu>1$, $\alpha_\nu^\ast\in[0,1]$ получаем
неравенство
\begin{equation}\label{ogrS}
0\leqslant S_\nu(\alpha)<1.
\end{equation}
Ясно что
\begin{equation}\label{ogrA}
0\leqslant A_{N+1}^{}(\alpha,t)<S_{N+1}^{}(\alpha).
\end{equation}
Используя (\ref{summa}) и (\ref{ogrS}) выводим оценку сверху
\begin{equation}\label{ogrG}
G_N(\alpha)<N.
\end{equation}
Из (\ref{s}) ,(\ref{ogrS}) и (\ref{ogrA}) получаем ограничение на
интеграл
\begin{equation}\label{ogrI}
 G_N(\alpha)\leqslant I_{\alpha}(t)<G_{N+1}(\alpha),
\end{equation}
и
\begin{equation}\label{DogrI}
I_{\alpha}(t)=G_N(\alpha)+O(1).
\end{equation}

Неравенство (\ref{ogrI}) дает оценку сверху $I_{\alpha}(t)<N+1$
для любого $\alpha\in\mathbb{R}$, из чего непосредственно следует
пункт 1) теоремы 2.

{\bf 3. Доказательство пункта 2) теоремы 2.}

Для доказательства пункта 2) теоремы 2 нам понадобятся непрерывные
дроби с {\it вещественными} неполными частными. В этом пункте и
далее для записи бесконечного множества аргументов будем
пользоваться обозначением $\overline{x}=(x_1,x_2,\ldots)$, где
$x_i\in[1;+\infty)$, $i\in\mathbb{N}$. Определим функцию
$\alpha(\overline{x})=[0;x_1,x_2,\ldots]$. Также определим функции
$q_\nu(\overline{x})$ и $p_\nu(\overline{x})$ как континуанты
$$
\begin{array}{cc}
  q_\nu(\overline{x})=\left\{%
\begin{array}{ll}
    \langle x_1,\ldots,x_\nu\rangle, & \hbox{$\nu\geqslant1;$} \\
    1, & \hbox{$\nu=0;$} \\
    0, & \hbox{$\nu=-1,$}
\end{array}%
\right. & p_\nu(\overline{x})=\left\{%
\begin{array}{ll}
    \langle x_1,\ldots,x_{\nu-1}\rangle, & \hbox{$\nu\geqslant1;$} \\
    0, & \hbox{$\nu=0;$} \\
    1, & \hbox{$\nu=-1.$}
\end{array}%
\right.\\
\end{array}
$$
Для всех $\nu\geqslant0$ выполняется (см \cite{Khinchin})
неравенство
\begin{equation}\label{star}
q_\nu(\overline{x})\geqslant2^{\frac{\nu-1}{2}}.
\end{equation}

Рассмотрим функцию
$\psi_{\alpha(\overline{x})}(t)=\|\alpha(\overline{x})
q_N^{}(\overline{x})\|$, интеграл
$I_{\alpha(\overline{x})}(t)=\int\limits_1^t\psi_{\alpha(\overline{x})}
(\xi)d\xi$ и функции
$$
    \alpha_\nu(\overline{x})=\left\{%
\begin{array}{ll}
    [x_\nu;x_{\nu+1},x_{\nu+2},\ldots], & \hbox{$\nu\geqslant1;$} \\
    \\

    [0;x_1,x_2,\ldots,], & \hbox{$\nu=0,$} \\
\end{array}%
\right.
\,\,\alpha_\nu^\ast(\overline{x})=\left\{%
\begin{array}{ll}
    [0;x_\nu,x_{\nu-1},x_{\nu-2},\ldots,x_1], & \hbox{$\nu\geqslant1;$} \\
    \\

    0, & \hbox{$\nu=0.$} \\
\end{array}%
\right.\\
$$

При каждом $\overline{x}$ выполняются следующие простейшие
свойства:
\begin{equation}\label{alpha}
\alpha_\nu(\overline{x})>1,   \,\, \nu\geqslant1,
\end{equation}

\begin{equation}\label{alast}
0\leqslant\alpha_\nu^\ast(\overline{x})\leqslant1,
\end{equation}

\begin{equation}\label{drob}
\alpha_\nu^\ast(\overline{x})=\frac{q_{\nu-1}(\overline{x})}{q_\nu(\overline{x})},
\end{equation}

\begin{equation}\label{alfn}
\alpha_\nu(\overline{x})=x_\nu+\frac{1}{\alpha_{\nu+1}(\overline{x})},
\end{equation}

\begin{equation}\label{alfnu}
\alpha_\nu^\ast(\overline{x})=\frac{1}{x_\nu+\alpha_{\nu-1}^\ast(\overline{x})},
\end{equation}

\begin{equation}\label{razn}
\alpha_\nu(\overline{x})-\frac{p_\nu(\overline{x})}{q_\nu(\overline{x})}=\frac{(-1)^\nu}
{q^2(\overline{x})(\alpha_{\nu+1}(\overline{x})+\alpha_\nu^\ast(\overline{x}))}.
\end{equation}

По аналогии с (\ref{ds}) и (\ref{summa}) рассмотрим функции
$$
S_\nu(\overline{x})=\frac{(1-\alpha_\nu^\ast(\overline{x}))\alpha_{\nu+1}(\overline{x})}
{\alpha_{\nu+1}(\overline{x})+\alpha_\nu^\ast(\overline{x})},\,\,G_n(\overline{x})=\sum\limits_{\nu=1}^n
S_\nu(\overline{x}).
$$

Аналогично формуле (\ref{DogrI}) можно получить равенство
\begin{equation}\label{DogrIotx}
I_{\alpha}(t)=G_N(\overline{x})+O(1),
\end{equation}
где $N$ определено аналогично (\ref{aa}).

Найдем частные производные от функций $\alpha_\nu(\overline{x})$ и
$\alpha_\nu^\ast(\overline{x})$ по переменной $x_k$. Для этого
уточним зависимость от $k$-ого аргумента:
$$
\alpha_{\nu}(\overline{x})=\left\{%
\begin{array}{ll}
    x_\nu+\frac{1}{\alpha_{\nu+1}(\overline{x})} & \hbox{$\nu\leqslant k-1;$} \\
    \\
  x_k+\frac{1}{\alpha_{\nu+1}(\overline{x})}, & \hbox{$\nu=k;$}\\
    \\
    x_\nu+\frac{1}{\alpha_{\nu+1}(\overline{x})}, & \hbox{$\nu\geqslant k+1,$} \\
\end{array}%
\right.
\alpha_{\nu}^\ast(\overline{x})=\left\{%
\begin{array}{ll}
    \frac{1}{x_\nu+\alpha_{\nu-1}^\ast(\overline{x})}, & \hbox{$\nu<k;$} \\
    \\
    \frac{1}{x_k+\alpha^\ast_{\nu-1}(\overline{x})},& \hbox{$\nu=k;$}\\
    \\
    \frac{1}{x_\nu+\alpha_{\nu-1}^\ast(\overline{x})}, & \hbox{$\nu\geqslant k+1.$} \\
\end{array}%
\right.
$$
Теперь для производных легко получить следующие соотношения:
\begin{equation}\label{pralpha}
\big(\alpha_{\nu}(\overline{x})\big)_{x_k}'=\left\{%
\begin{array}{ll}
    -\frac{\big(\alpha_{\nu+1}(\overline{x})\big)_{x_k}'}{\big(\alpha_{\nu+1}(\overline{x})\big)^2},
    & \hbox{$\nu\leqslant k-1;$} \\
    \\
  1, & \hbox{$\nu=k;$}\\
    \\
    0, & \hbox{$\nu\geqslant k+1,$} \\
\end{array}%
\right.
\end{equation}

\begin{equation}\label{pralfast}
\big(\alpha^\ast_{\nu}(\overline{x})\big)_{x_k}'=\left\{%
\begin{array}{ll}
    0, & \hbox{$\nu<k;$} \\
    \\
    -\frac{1}{\big(x_k+\alpha^\ast_{\nu-1}(\overline{x})\big)^2},& \hbox{$\nu=k;$}\\
    \\
    -\frac{\big(\alpha_{\nu-1}^\ast(\overline{x})\big)_{x_k}'}{\big(x_{\nu}+\alpha_{\nu-1}^\ast(\overline{x})\big)^2},
     & \hbox{$\nu\geqslant k+1.$} \\
\end{array}%
\right.
\end{equation}

Пусть $l\in\mathbb{N}_0$. Из (\ref{pralpha}) получаем
\begin{equation}\label{znpra}
{\rm sign}\left(\alpha_{\nu}(\overline{x})\right)_{x_k}'=\left\{%
\begin{array}{ll}
    1, & \hbox{$\nu=k-2l;$} \\
    \\
    -1, & \hbox{$\nu=k-1-2l.$}\\
\end{array}%
\right. \end{equation}

Аналогично и для $(\alpha^\ast_{\nu}(\overline{x}))_{x_k}'$ из
(\ref{pralfast}) видно, что
\begin{equation}\label{znprast}
{\rm sign}\left(\alpha_{\nu}^\ast(\overline{x})\right)_{x_k}'=\left\{%
\begin{array}{ll}
    -1, & \hbox{$\nu=k+2l;$} \\
    \\
    1, & \hbox{$\nu=k+1+2l.$}\\
\end{array}%
\right.
\end{equation}

{\bf Лемма 1.} {\it Функция $G_n(\overline{x})$ возрастает по
каждому из первых $n+1$ аргументу}.

{\bf Доказательство}: Покажем, что
$\big(G_n(\overline{x})\big)_{x_k}'>0$, для $k\leqslant n+1$.
Найдем производную по $k$-ому аргументу. В этом доказательстве не
будем писать аргумент $\overline{x}$.
\begin{equation}\label{astr}
\big(G_n(\overline{x})\big)_{x_k}'=\sum\limits_{\nu=1}^n\frac{(\alpha_{\nu+1})_{x_k}'(1-\alpha_{\nu}^\ast)
\alpha_{\nu}^\ast}{(\alpha_{\nu+1}+\alpha_{\nu}^\ast)^2}
+\sum\limits_{\nu=1}^n\frac{(-\alpha_{\nu}^\ast)_{x_k}'\alpha_{\nu+1}(\alpha_{\nu+1}+1)}
{(\alpha_{\nu+1}+\alpha_{\nu}^\ast)^2}.
\end{equation}

Исследуем первую сумму из (\ref{astr}) и покажем, что она
неотрицательная.

Воспользовавшись соотношением (\ref{pralpha}), отбросим нулевые
слагаемые
$$
\sum\limits_{\nu=1}^n\frac{(\alpha_{\nu+1})_{x_k}'(1-\alpha_\nu^\ast)\alpha_\nu^\ast}{(\alpha_{\nu+1}+\alpha_\nu^\ast)^2}=
\sum\limits_{\nu=1}^{k-1}\frac{(\alpha_{\nu+1})_{x_k}'(1-\alpha_\nu^\ast)\alpha_\nu^\ast}
{(\alpha_{\nu+1}+\alpha_\nu^\ast)^2}.
$$
В новой сумме, используя условие $\alpha_0^\ast=0$, можно при
необходимости сделать четное количество слагаемых, добавив
слагаемое с $\nu=0$. Сгруппируем слагаемые по два, начиная с
последнего, и преобразуем сумму, воспользовавшись (\ref{alfn}) и
(\ref{pralpha}):
$$
\sum\limits_{\nu=1}^{k-1}\frac{(\alpha_{\nu+1})_{x_k}'(1-\alpha_\nu^\ast)\alpha_\nu^\ast}
{(\alpha_{\nu+1}+\alpha_\nu^\ast)^2}=\sum\limits_{l=0}^{\left[\frac{k-2}{2}\right]}\left(\frac{(1-\alpha^\ast_{k-1-2l})
\alpha^\ast_{k-1-2l}(\alpha_{k-2l})_{x_k}'}
{(\alpha_{k-2l}+\alpha^\ast_{k-1-2l})^2}
+\frac{(1-\alpha^\ast_{k-2-2l})\alpha^\ast_{k-2-2l}(\alpha_{k-1-2l})_{x_k}'}{(\alpha_{k-1-2l}+\alpha^\ast_{k-2-2l})^2}
\right)=
$$
$$
=\sum\limits_{l=0}^{\left[\frac{k-2}{2}\right]}\left(\frac{\left(1-\frac{1}{x_{k-1-2l}+\alpha_{k-2-2l}^\ast}\right)
\frac{1}{x_{k-1-2l}+\alpha_{k-2-2l}^\ast}}
{\left(\alpha_{k-2l}+\frac{1}{x_{k-1-2l}+\alpha_{k-2-2l}^\ast}\right)^2}(\alpha_{k-2l})_{x_k}'-
\frac{(1-\alpha^\ast_{k-2-2l}) \alpha^\ast_{k-2-2l}}
{\left(x_{k-1-2l}+\frac{1}{\alpha_{k-2l}}+\alpha^\ast_{k-2-2l}\right)^2}\frac{(\alpha_{k-2l})_{x_k}'}
{(\alpha_{k-2l})^2}\right)=
$$
$$
=\sum\limits_{l=0}^{\left[\frac{k-2}{2}\right]}(\alpha_{k-2l})_{x_k}'
\left(\frac{x_{k-1-2l}-1+(\alpha_{k-2-2l}^\ast)^2}{(\alpha_{k-2l}(x_{k-1-2l}+
\alpha_{k-2-2l}^\ast)+1)^2}\right).
$$

Из (\ref{znpra}), следует что $(\alpha_{k-2l})_{x_k}'>0$, для всех
$l$, а значит и вся сумма тоже положительна.

Покажем, что вторая сумма из (\ref{astr}) всегда не отрицательна.

Из (\ref{pralfast}) и (\ref{znprast}) получаем, что сумма
знакочередующаяся и начинается с положительного слагаемого при
$\nu=k$. Сгруппируем слагаемые по два. Если в сумме нечетное
количество слагаемых, то последнее слагаемое можно отбросить, ибо
оно положительно, от этого сумма только уменьшится. Воспользуемся
(\ref{alfn}), (\ref{pralfast}), (\ref{znprast}) и оценим сумму
снизу
$$
\sum\limits_{\nu=1}^n\frac{(1+\alpha_{\nu+1})\alpha_{\nu+1}(-\alpha^\ast_\nu)_{x_k}'}
{(\alpha_{\nu+1}+\alpha^\ast_\nu)^2}=\sum\limits_{\nu=k}^n\frac{(1+\alpha_{\nu+1})\alpha_{\nu+1}(-\alpha^\ast_\nu)_{x_k}'}
{(\alpha_{\nu+1}+\alpha^\ast_\nu)^2}\geqslant
$$
$$
\geqslant\sum\limits_{l=0}^{\left[\frac{n-k-1}{2}\right]}
\left(-\frac{(1+\alpha_{k+1+2l})\alpha_{k+1+2l}
(\alpha^\ast_{k+2l})_{x_k}'}{(\alpha_{k+1+2l}+\alpha^\ast_{k+2l})^2}-\frac{(1+\alpha_{k+2+2l})\alpha_{k+2+2l}
(\alpha^\ast_{k+1+2l})_{x_k}'}{(\alpha_{k+2+2l}+\alpha^\ast_{k+1+2l})^2}\right)=
$$
$$
\sum\limits_{l=0}^{\left[\frac{n-k-1}{2}\right]}
\left(-\frac{\left(1+x_{k+1+2l}+\frac{1}{\alpha_{k+2+2l}}\right)\left(x_{k+1+2l}+\frac{1}{\alpha_{k+2+2l}}\right)
(\alpha^\ast_{k+2l})_{x_k}'}{\left(x_{k+1+2l}+
\frac{1}{\alpha_{k+2+2l}}+\alpha^\ast_{k+2l}\right)^2}-\right.
$$
$$
\left.-\frac{(1+\alpha_{k+2+2l})\alpha_{k+2+2l}}
{\left(\alpha_{k+2+2l}+\frac{1}
{x_{k+1+2l}+\alpha_{k+2l}^\ast}\right)^2}\frac{(-\alpha_{k+2l}^\ast)_{x_k}'}{(x_{k+1+2l}+\alpha_{k+2l}^\ast)^2}\right)=
$$
$$
=\sum\limits_{l=0}^{\left[\frac{n-k-1}{2}\right]}
(-\alpha^\ast_{k+2l})_{x_k}'\left(\frac{(\alpha_{k+2+2l}+x_{k+1+2l}\alpha_{k+2+2l}+1)(x_{k+1+2l}\alpha_{k+2+2l}+1)
-(1+\alpha_{k+2+2l})
\alpha_{k+2+2l}}{(\alpha_{k+2+2l}(x_{k+1+2l}+\alpha_{k+2l}^\ast)+1)^2}\right)>0.
$$
Получили, что производная положительна при $k\leqslant n+1$,
значит функция $G_n(\overline{x})$ возрастает по первым $n+1$
аргументам. Доказательство леммы 1 завершено.
\\
Из леммы 1 сразу получаем

{\bf Следствие 1.} {\it Для любого набора $\overline{x}$ и
$n\in\mathbb{N}$ выполняется неравенство}
\begin{equation}\label{Gapr}
G_n(\overline{x})\geqslant
G_n(\underbrace{1,1,\ldots,1}\limits_{n+1}\,,x_{n+2},\ldots).
\end{equation}

Возьмем вещественное число $z$, для него рассмотрим набор
$\overline{z}=(z,z,\ldots)$ и определим число
\begin{equation}\label{sz}
S(z)=\frac{(1-\alpha(\overline{z}))(z+\alpha(\overline{z}))}{z+2\alpha(\overline{z})}.
\end{equation}

{\bf Лемма 2.} {\it Для любого $z\in[1,+\infty)$ выполняется
неравенство}
$$
|G_n(\overline{z})-S(z)n|\leqslant 4.
$$

{\bf Доказательство:} Для $\overline{z}=(z,z,\ldots)$ выполняется
свойство
$\alpha_\nu^\ast(\overline{z})=\frac{p_\nu(\overline{z})}{q_\nu(\overline{z})}$,
которые позволяют записать разность (\ref{razn}) следующим образом
$$
\alpha(\overline{z})-\alpha_\nu^\ast(\overline{z})=
\frac{(-1)^\nu}{q^2_\nu(\overline{z})(\alpha_{\nu+1}(\overline{z})+\alpha_\nu^\ast(\overline{z}))}.
$$
Запишем $G_n(\overline{z})$ следующим образом
$$
G_n(\overline{z})=\sum\limits_{\nu=1}^n
S(z)+\sum\limits_{\nu=1}^n(S_\nu(\overline{z})-S(z)).
$$
Воспользуемся (\ref{star}) и покажем, что ряд
$\sum\limits_{\nu=1}^\infty(S_\nu(\overline{x})-S(z))$ сходится
абсолютно
$$
\sum\limits_{\nu=1}^\infty|S_\nu(\overline{z})-S(z)|=\sum\limits_{\nu=1}^\infty\left|\frac{(1-\alpha_\nu^\ast(\overline{z}))
\alpha_{\nu+1}(\overline{z})}{\alpha_{\nu+1}(\overline{z})+\alpha_\nu^\ast(\overline{z})}-\frac{(1-\alpha(\overline{z}))
\alpha_{\nu+1}(\overline{z}))}{\alpha_{\nu+1}(\overline{z})+\alpha(\overline{z})}\right|=
$$
$$
=\sum\limits_{\nu=1}^\infty\frac{\alpha_{\nu+1}(\alpha_{\nu+1}+1)}{(\alpha_{\nu+1}+\alpha_\nu^\ast)^2
(\alpha_{\nu+1}+\alpha)q_\nu^2(\overline{z})}\leqslant\sum\limits_{\nu=1}^\infty\frac{2}{q_\nu^2(\overline{z})}\leqslant
\sum\limits_{\nu=1}^\infty\frac{4}{2^\nu}=4,
$$
откуда
$$
|G_n(\overline{z})-S(z)n|\leqslant4,
$$ что и требовалось доказать.

{\bf Лемма 3.} {\it Пусть
$\overline{x}=(z_1,z_2,\ldots,z_{n+1},x_{n+2},\ldots)$ и
$\overline{y}=(z_1,z_2,\ldots,z_{n+1},y_{n+2},\ldots)$, тогда}
$$
|G_n(\overline{x})-G_n(\overline{y})|<1.
$$

{\bf Доказательство.} При $\nu\leqslant n+1$ для функций
$\alpha_\nu^\ast(\overline{x})$ и $\alpha_\nu^\ast(\overline{y})$
выполняется равенство

$$
\alpha_\nu^\ast(\overline{x})=\alpha_\nu^\ast(\overline{y}).
$$

Так как
$\alpha_\nu(\overline{x})=[x_\nu;x_{\nu+1},\ldots,x_{n+1},a_{n+2},\ldots]$
и
$\alpha_\nu(\overline{y})=[x_\nu;x_{\nu+1},\ldots,x_{n+1},b_{n+2},\ldots]$,
то значения обеих функций лежат между числами
$\frac{p_{n-\nu+1}(\overline{x})}{q_{n-\nu+1}(\overline{x})}$ и
$\frac{p_{n-\nu+1}(\overline{x})+p_{n-\nu} (\overline{x})}
{q_{n-\nu+1}(\overline{x})+q_{n-\nu}(\overline{x})}$ (см. \cite
{Khinchin}). Неравенство (\ref{star}) позволяет оценить расстояние между ними:
$$
|\alpha_\nu(\overline{x})-\alpha_\nu(\overline{y})|<\frac{1}{q_{n-\nu+1}(\overline{x})(q_{n-\nu}(\overline{x})
+q_{n-\nu+1}(\overline{x}))}\leqslant\frac{1}{2^{n-\nu}}.
$$
Оценим разность
$$
|S_\nu(\overline{x})-S_\nu(\overline{y})|=\left|\frac{(1-\alpha_\nu^\ast(\overline{x}))\alpha_{\nu+1}(\overline{x})}
{\alpha_\nu^\ast(\overline{x})+\alpha_{\nu+1}(\overline{x})}-\frac{(1-\alpha_\nu^\ast(\overline{y}))\alpha_{\nu+1}
(\overline{y})}{\alpha_\nu^\ast(\overline{y})+\alpha_{\nu+1}(\overline{y})}\right|=
$$
$$=
(1-\alpha_\nu^\ast(\overline{x}))\left|\frac{\alpha_{\nu+1}(\overline{x})\alpha_\nu^\ast(\overline{x})-
\alpha_{\nu+1}(\overline{y})\alpha_\nu^\ast(\overline{x})}{(\alpha_\nu^\ast(\overline{x})+\alpha_{\nu+1}(\overline{x}))
(\alpha_\nu^\ast(\overline{x})+\alpha_{\nu+1}(\overline{y}))}\right|\leqslant
$$
$$
\leqslant(1-\alpha_\nu^\ast(\overline{x}))\alpha_\nu^\ast(\overline{x})|\alpha_{\nu+1}(\overline{x})-
\alpha_{\nu+1}(\overline{y})|
\leqslant\frac{1}{4}|\alpha_{\nu+1}(\overline{x})-\alpha_{\nu+1}(\overline{y})|\leqslant\frac{1}{2^{n-\nu+1}}.
$$

Далее
$$
|G_n(\overline{z})-G_n(\overline{y})|=\left|\sum\limits_{\nu=1}^nS_\nu(\overline{z})-\sum\limits_{\nu=1}^nS_\nu
(\overline{y})\right|
\leqslant\sum\limits_{\nu=1}^n|S_\nu(\overline{z})-S_\nu(\overline{y})|
\leqslant\sum\limits_{\nu=1}^n\frac{1}{2^{n-\nu+1}}<1.
$$
Доказательство леммы 3 завершено.

{\bf Следствие 2.} {\it Для
$\overline{y}=(\underbrace{z,\ldots,z}_{n+1},y_{n+2}\ldots)$ из
лемм 2 и 3 получаем что}
\begin{equation}\label{GnSz}
|G_n(\overline{y})-S(z)n|<5.
\end{equation}

Теперь мы завершим доказательство пункта 2) теоремы 2.
Из(\ref{DogrI}), (\ref{Gapr}), (\ref{sz}) и (\ref{GnSz}) получаем,
что для любого набора $\overline{x}$ выполняется
$$
\liminf\limits_{t\rightarrow\infty}\frac{I_{\alpha(\overline{x})}(t)}{N}=
\liminf\limits_{t\rightarrow\infty}\frac{G_N^{}(\overline{x})}{N}=
\liminf\limits_{n\rightarrow\infty}\frac{G_n(\overline{x})}{n}\geqslant
\liminf\limits_{n\rightarrow\infty}\frac{G_n(1,1,\ldots,1,y_{n+2},\ldots)}{n}=S(1)=\left(\frac{1}{2}-
\frac{\sqrt{5}}{10}\right).
$$

Доказательство теоремы 2 завершено.

{\bf 4. Доказательство  теоремы 3.}

{\bf Лемма 4}. Пусть $\overline{x}=(0;x,z_2,z_3\ldots,)$,
$\overline{y}=(0;y,z_2,z_3\ldots,)$ и $n\in\mathbb{N}$ тогда
$$
\left|G_n(\overline{x})-G_n(\overline{y})\right|<8.
$$

{\bf Доказательство.} Оба числа $\alpha_\nu^\ast(\overline{x})$ и
$\alpha_\nu^\ast(\overline{y})$ попадают в интервал между числами
$\frac{p_{\nu-1}(\overline{x})}{q_{\nu-1}(\overline{x})}$ и
$\frac{p_{\nu-1}(\overline{x})+p_{\nu-2}(\overline{x})}
{q_{\nu-1}(\overline{x})+q_{\nu-2}(\overline{x})}$ (см. \cite
{Khinchin}), где $\frac{p_\nu(\overline{x})}{q_\nu(\overline{x})}$
- подходящая дробь к $\alpha_\nu^\ast(\overline{x})$. Это
позволяет оценить расстояние между ними:
$$
|\alpha_\nu^\ast(\overline{x})-\alpha_\nu^\ast(\overline{y})|<\frac{1}{q_{\nu-1}(\overline{x})(q_{\nu-2}(\overline{x})
+q_{n-1}(\overline{x}))}\leqslant\frac{1}{2^{\nu-2}}.
$$

При $\nu\geqslant 2$ выполняется равенство
$$
\alpha_\nu(\overline{x})=\alpha_\nu(\overline{y}).
$$

Оценим разность
$$
|G_n(\overline{x})-G_n(\overline{y})|=\left|\sum\limits_{\nu=1}^nS_\nu(\overline{x})-\sum\limits_{\nu=1}^nS_\nu
(\overline{y})\right|
\leqslant\sum\limits_{\nu=1}^n|S_\nu(\overline{x})-S_\nu(\overline{y})|=
$$
$$
=\sum\limits_{\nu=1}^n\left|\frac{(1-\alpha_\nu^\ast(\overline{x}))\alpha_{\nu+1}(\overline{x})}
{\alpha_\nu^\ast(\overline{x})+\alpha_{\nu+1}(\overline{x})}-\frac{(1-\alpha_\nu^\ast(\overline{y}))
\alpha_{\nu+1}(\overline{y})}{\alpha_\nu^\ast(\overline{y})+\alpha_{\nu+1}(\overline{y})}\right|=
$$
$$
=\sum\limits_{\nu=1}^n\alpha_{\nu+1}(\overline{x})\left|\frac{\alpha_{\nu}^\ast(\overline{y})-
\alpha_\nu^\ast(\overline{x})\alpha_{\nu+1}(\overline{y})-\alpha_\nu^\ast(\overline{x})+\alpha_\nu^\ast(\overline{y})
\alpha_{\nu+1}(\overline{x})}{(\alpha_\nu^\ast(\overline{x})+\alpha_{\nu+1}(\overline{x}))
(\alpha_\nu^\ast(\overline{y})+\alpha_{\nu+1}(\overline{y}))}\right|\leqslant
$$
$$
\leqslant\sum\limits_{\nu=1}^n\frac{\alpha_{\nu+1}(\overline{x})(1+\alpha_{\nu+1}(\overline{y}))}
{\alpha_{\nu+1}(\overline{x})\alpha_{\nu+1}(\overline{y})}\left|\alpha_{\nu}^\ast(\overline{y})-
\alpha_\nu^\ast(\overline{x})\right|\leqslant2\sum\limits_{\nu=1}^n\left|\alpha_{\nu}^\ast(\overline{y})-
\alpha_\nu^\ast(\overline{x})\right|\leqslant\sum\limits_{\nu=1}^\infty\frac{8}{2^\nu}=8.
$$
Доказательство леммы 4 завершено.

Из следствия 2 и леммы 4 получаем, что для набора
$\overline{x}=(x_1,\underbrace{z,z,\ldots,z}_{n},x_{n+2},\ldots)$
выполняется неравенство
\begin{equation}\label{odnzam}
\left|G_n(\alpha)-S(z)n\right|<13.
\end{equation}

Сумма $G_n(\alpha)$ зависит от $n+1$ аргумента $(a_1,a_2,\ldots,
a_n,\alpha_{n+1})$. В дальнейшем рассуждении нам удобно выделить
зависимость от последнего аргумента. Для числа
$\alpha=[0;a_1,a_2,\ldots]$ через $G_n(\alpha,x)$ будем обозначать
сумму для набора $(a_1,a_2,\ldots,a_n,x)$, где
$x\in(1,+\infty)$. Обозначим
$G_n(\alpha,+\infty)=\lim\limits_{x\rightarrow\infty}G_n(\alpha,x)$
и $G_n(\alpha,1)=\lim\limits_{x\rightarrow1}G_n(\alpha,x)$.

{\bf Утверждение.} Для любых $x,y\in(1,+\infty)$ выполняется неравенство
\begin{equation}\label{dobav}
|G_{n-1}(\alpha,x)-G_n(\alpha,y)|<3.
\end{equation}

{\bf Доказательство.} Из леммы 1 и леммы 3 получаем
$$0<G_n(\alpha,+\infty)-G_n(\alpha,1)<1.$$

Из формулы (\ref{summa}) можно вывести равенство
$$G_n(\alpha,+\infty)=G_{n-1}(\alpha,a_n)+(1-\alpha_n^\ast).$$

Объединив эти формулы, получаем неравенство
$$|G_{n-1}(\alpha,a_n)-G_n(\alpha,1)|<1.$$

Применим лемму 3 для каждой из сумм и получим утверждение.

Для функции $S(z)$ выполняются свойства:

1) она монотонно возрастает;

2) $\frac{1}{2}-\frac{\sqrt{5}}{10}=S(1)\leqslant S(z)\leqslant
S(+\infty)=1.$

Для любого $d=\left[\frac{1}{2}-\frac{\sqrt{5}}{10},1\right)$
можно взять натуральные числа $a$ и $b$, такие чтобы выполнялось
условие $S(a)\leqslant d<S(b)$.

{\bf Лемма 5.} Существует $t_{min}$ такое, что для любого набора
$(a_1,\ldots,a_k)$ и любого $x\in[1,+\infty)$ при всех
$t>t_{min}$ выполняются неравенства:

1) $$
\frac{G_{k+t}(\alpha,x)}{k+t}<d,
$$
для всех чисел вида
$\alpha=[a_1,\ldots,a_k,\underbrace{a,\ldots,a}_t,x];$

2) $$ \frac{G_{k+t}(\beta,x)}{k+t}>d,
$$
для всех чисел вида
$\beta=[a_1,\ldots,a_k,\underbrace{b,\ldots,b}_t,x]$.

{\bf Доказательство.} Запишем сумму $G_{k+t}$ следующим образом
$$
G_{k+t}=G_k+\sum_{\nu=k+1}^{k+t}S_\nu,
$$
Сумма $\sum\limits_{\nu=k+1}^{k+t}S_\nu$ есть ни что иное как
$G_t$ для числа
$(\frac{1}{\alpha_{n+1}^\ast},\underbrace{a,\ldots,a}_{t-1},x)$.
Запишем $G_{k+t}$ через две суммы, а потом применим лемму 4 и
получим равенство
\begin{equation}\label{formula0}
G_{k+t}(\alpha,x)=G_{k}+G_{t}\left(\frac{1}{\alpha_{n+1}^\ast},\underbrace{a,\ldots,a}_{t-1},x\right)=G_k+S(a)(t-1)+R,
\end{equation}
где $R<13$. Рассмотрим предел
$$
\lim_{t\rightarrow+\infty}\frac{G_{k+t}}{k+t}=\lim_{t\rightarrow+\infty}\frac{G_k+S(a)(t-1)+R}{k+t}=S(a)<d.
$$
Откуда следует утверждение леммы. Пункт 2 доказывается аналогично.

Для последовательности натуральных чисел $n_\nu$ будем обозначать
частичную сумму через $W_t=\sum\limits_{\nu=1}^tn_\nu$.

{\bf Следствие 3.} Существует
$\alpha=[0;\underbrace{a,\ldots,a}_{n_1},\underbrace{b,\ldots,b}_{n_2},\underbrace{a,\ldots,a}_{n_3},\ldots]$,
такое что
$$
\frac{G_{n_1}(\alpha,+\infty)}{n_1}<d,
$$
и для четного $t$ выполняется
$$
\frac{G_{W_t-1}(\alpha,1)}{W_t-1}<d<\frac{G_{W_t}(\alpha,1)}{W_t},
$$
а для нечетного $t>1$ выполняется
$$
\,\,\frac{G_{W_t}(\alpha,+\infty)}{W_t}<d<\frac{G_{W_t-1}(\alpha,+\infty)}{W_t-1}.
$$
Замечание. При доказательстве следствия 3 в качестве $n_\nu$ надо
брать $t_{min}$ из леммы 5.

{\bf Доказательство теоремы 3.} Возьмем число $\alpha$ из
следствия 3 и покажем, что все $n_{t+1}$ ограничены. Пусть $t$
нечетно, тогда
$$
\frac{G_{W_t}(\alpha,1)}{W_t}<\frac{G_{W_t}(\alpha,+\infty)}{W_t}<d,
$$
$$
\frac{G_{W_{t+1}-1}(\alpha,1)}{W_{t+1}-1}<d.
$$
Обозначим $\tilde{b}=[\underbrace{b;b,\ldots,b}_{n_{t+1}-1},1]$.
По лемме 3 получаем
$$G_{W_t}(\alpha,\tilde{b})=G_{W_t}(\alpha,1)+R_1,$$ где $|R_1|<1.$
Воспользуемся (\ref{formula0}) и распишем
$$
G_{W_{t+1}-1}(\alpha,1)=G_{W_t}(\alpha,\tilde{b})+S(b)(n_{t+1}-1)+R=G_{W_t}(\alpha,1)+S(b)n_{t+1}+R+R_1-S(b),
$$

$$
\frac{G_{W_t}(\alpha,1)+S(b)n_{t+1}+R+R_1-S(b)}{W_t+n_{t+1}-1}=\frac{G_{W_{t+1}-1}(\alpha,1)}{W_{t+1}-1}<d,
$$
и получаем
$$
n_{t+1}<\frac{dW_t-G_{W_t}(\alpha,1)-d-R-R_1+S(b)}{S(b)-d}.
$$

Пусть $t$ четно, тогда
$$
d<\frac{G_{W_t}(\alpha,1)}{W_t}<\frac{G_{W_t}(\alpha,+\infty)}{W_t},
$$
$$
d<\frac{G_{W_{t+1}-1}(\alpha,+\infty)}{W_{t+1}-1}.
$$
Обозначим
$\tilde{a}=\lim\limits_{x\rightarrow+\infty}[\underbrace{a;a,\ldots,a}_{n_{t+1}-1},x]=
[\underbrace{a;a,\ldots,a}_{n_{t+1}-1}]$.
По лемме 3 и (9) для $\tilde{a}$ получаем
$$G_{W_t}(\tilde{a})=G_{W_t}(\alpha,1)+R_2,$$ где $|R_2|<1.$
Подставим
$$
G_{W_{t+1}-1}(\alpha,+\infty)=G_{W_t}(\alpha,\tilde{a})+S(a)(n_{t+1}-1)+R=G_{W_t}(\alpha,1)+S(a)n_{t+1}+R+R_2-S(a),
$$
$$
d<\frac{G_{W_{t+1}-1}(\alpha,+\infty)}{W_{t+1}-1}=\frac{G_{W_t}(\alpha,1)+S(a)n_{t+1}+R+R_2-S(a)}{W_t+n_{t+1}-1},
$$
$$
n_{t+1}<\frac{d+R+R_2+S(a)+G_{W_t}(1)-dW_t}{d-S(a)}.
$$
В обоих случаях все $n_{t+1}$ ограничены некоторой константой,
назовем её $M$.

Для $\frac{G_{n}(\alpha)}{n}$ и $\frac{G_{n+k}(\alpha)}{n+k}$ из
(\ref{summa}), (\ref{ogrS}) и (\ref{ogrG}) выводим неравенство
\begin{equation}\label{ogrRazG}
\left|\frac{G_{n+k}}{n+k}-\frac{G_{n}}{n}\right|=\left|\frac{G_{n}+\sum\limits_{\nu=n+1}^{n+k}S_\nu}
{n+k}-\frac{G_{n}}{n}\right|=\left|\frac{\sum\limits_{\nu=n+1}^{n+k}S_\nu-k\frac{G_{n}}{n}}
{n+k}\right|<\frac{k}{n+k}.
\end{equation}

Для последовательности $\frac{G_n(\alpha)}{n}$ выполняются
свойства:

1) $\left|\frac{G_n}{n}-\frac{G_{n+1}}{n+1}\right|<\frac{1}{n}$ из (\ref{ogrRazG}),

2) $\left|\frac{G_{W_t}}{W_t}-d\right|<\frac{3}{W_t}$ из утверждения и следствия 3,

3) $W_{t+1}-W_t<M$ из ограниченности $n_t$.

Пусть $n\in[W_t;W_{t+1})$, тогда из свойств 1 и 3
получаем $\left|\frac{G_n}{n}-d\right|<\frac{M+3}{W_t}$. Для $n>M$
получим, что $n-M<W_t$ и
$\left|\frac{G_n}{n}-d\right|<\frac{M+3}{n-M}$.

Возьмем предел от правой и левой частей
$$
\lim_{n\rightarrow+\infty}\left|\frac{G_n}{n}-d\right|<\lim_{n\rightarrow+\infty}\frac{M+3}{n-M}=0,
$$
или
$$
\lim_{n\rightarrow+\infty}\frac{G_n}{n}=d.
$$

Случай $d=1$ получается для числа, у которого неполные частные
образуют возрастающую последовательность. Теорема 3 доказана.

{\bf 5. Эргодические свойства преобразования Гаусса.}

Нам понадобятся некоторые сведения из эргодической теории.
Основные нужные нам понятия и утверждения имеются в
\cite{Kornfeld} (см. также \cite{nakada,nakada1}).

Рассмотрим преобразование Гаусса  $T:[0;1)\rightarrow[0;1)$,
которое есть эндоморфизм задаваемый формулой
$$
Tx=\left\{%
\begin{array}{ll}
    \left\{\frac{1}{x}\right\}, & \hbox{при $x\neq0$}; \\
\\
    0, & \hbox{при $x=0$.} \\
\end{array}%
\right.
$$
Если для $x$ известно его разложение в цепную дробь
$x=[0;a_1,a_2,a_3,\ldots]$, то
$$
Tx=[0,a_2,a_3,\ldots].
$$
Инвариантная мера для преобразования Гаусса задается формулой
$$
\mu(A)=\frac{1}{\ln2}\int\limits_A\frac{1}{1+x}dx.
$$

Естественным расширением для преобразования Гаусса будет
автоморфизм $\widehat{T}:[0;1)^2\rightarrow[0;1)^2$, определяемый
как
$$
 \widehat{T}(x,y)=\left\{%
\begin{array}{ll}
   \left(\left\{\frac{1}{x}\right\},\frac{1}{\left[\frac{1}{x}\right]+y}\right),& \hbox{ при $x\neq0$;}\\
   \\
    (0,y), & \hbox{при $x=0$.} \\
\end{array}%
\right.
$$
У естественного расширения $\widehat{T}$ инвариантная мера есть
$$
\mu_2(A)=\frac{1}{\ln2}\int\int\limits_A\frac{1}{(1+xy)^2}dxdy.
$$

Преобразование $\widehat{T}$ обладает свойством $K$-перемешивания.
В  частности, преобразование\\ $\widehat{T}$ эргодично, и, согласно
теореме Биркгофа-Хинчина, для любой абсолютно интегрируемой
функции $f(x,y)$  асимптотическое равенство
\begin{equation}\label{int}
\lim_{n\to +\infty} \sum_{\nu =1}^n f(\widehat{T}^\nu(x,y)) =
\frac{1}{\ln 2} \int_0^1\int_0^1 \frac{f(x,y)dxdy}{(1+xy)^2}
\end{equation}
будет выполнено для почти всех $(x,y) \in [0,1)^2$.

Если на точку $(x,y)$ подействовать преобразованием
$\widehat{T}^\nu$, то
\begin{equation}\label{Wxy}
\widehat{T}^\nu(x,y)=\left(T^\nu(x),\frac{q_{\nu-1}+yp_{\nu-1}}{q_\nu+yp_\nu}\right).
\end{equation}
где $\frac{p_\nu}{q_\nu}$ -
 подходящие дроби для $x$.

Рассмотрим преобразование $\widehat{\widehat{T{}
}}:[0,1)^2\rightarrow[0,1)^2$, которое определим так
$$
\widehat{\widehat{T{}}}\left(\begin{array}{c}
                         z_1 \\
                         z_2\\
                       \end{array}\right)=\left(\begin{array}{c}
                         \widehat{T}(z_1) \\
                         \widehat{T}(z_2)\\
                       \end{array}\right), \, z_1,z_2\in[0,1)^2.
$$
Его инвариантная мера есть $\mu_2(z_1)\times\mu_2(z_2)$.
Преобразование $\widehat{\widehat{T{}}}$ эргодично, это следует из
того, что $\widehat{T}$ обладает свойством перемешивания (см.
\cite{Furs}).

Следующая теорема доказана Халасом (см. \cite{Hal}).

{\bf Теорема Халаса.}\,\,{\it Для любой интегрируемой функции
$\varphi(p)$ и любого эргодического\\
преобразования $T$
пространства $R$ конечной меры выражение
$$
\sum_{\nu =0}^{n-1} \varphi(T^{\nu}
p)-n\int\limits_R\varphi(p)d\mu
$$
для почти всех $p$, меняет знак бесконечное число раз в слабом
смысле, т.е. эта разность не может быть постоянно положительной
или отрицательной. }

Применяя эту теорему к эргодическому преобразованию
$\widehat{\widehat{T{}}}$, получим

{\bf Следствие 3.}\,\,{\it  Для почти всех
$(z_1,z_2)=(x_1,y_1,x_2,y_2)\in[0,1)^2\times[0,1)^2$ и для любой
$\mu_2$- интегрируемой функции $f(z)=f(x,y)$ разностная функция
$$
\sum\limits_{\nu=1}^nf(\widehat{T}^\nu
z_1)-\sum\limits_{\nu=1}^nf(\widehat{T}^\nu z_2),
$$
меняет знак бесконечное число раз в слабом смысле.}

{\bf 6. Доказательства Теоремы 1.} Рассмотрим функцию
\begin{equation}\label{fufu}
f(x,y)=\frac{1-y}{1+xy}. \end{equation} Поскольку
$$
\int\limits_0^1\int\limits_0^1\frac{1-y}{(1+xy)^3}dxdy=\frac{\ln2}{2},
$$
то применяя теорему Биркгофа-Хинчина (\ref{int}), приходим к
равенству
\begin{equation}\label{osn}
\lim\limits_{n\rightarrow\infty}\frac{1}{n}\sum\limits_{\nu=1}^n
f(\widehat{T}^\nu(x,y))=\frac12.
\end{equation}

Отметим, что если $(x,y)=\left(\frac{1}{\alpha_1},0\right)$, то
$f(\widehat{T}^\nu(x,0))=\frac{(1-\alpha_\nu^\ast)}{1+\frac{\alpha_\nu^\ast}{\alpha_{\nu+1}}}$,
и  $G_n(\alpha)=\sum\limits_{\nu=1}^n f(\widehat{T}^\nu(x,0))$.
Следующее утверждение, близко к использовавшемуся в работе
\cite{Kraaim}

{\bf Лемма 5.} {\it Для любого $(x,y)\in[0;1)^2$ и любого
$n\in\mathbb{N}$ выполняется неравенство}
$$
\sum\limits_{\nu=1}^n \left|f(\widehat{T}^\nu(x,y))-
f(\widehat{T}^\nu(x,0))\right|<4.
$$

{\bf Доказательство.} Покажем, что ряд
$$
\sum\limits_{\nu=1}^\infty
\left(f(\widehat{T}^\nu(x,0))-f(\widehat{T}^\nu(x,y))\right)
$$
сходится абсолютно. Обозначим
$\tilde{\alpha}_\nu=\frac{q_{\nu-1}+yp_{\nu-1}}{q_\nu+yp_\nu}$.
Далее
$$\sum\limits_{\nu=1}^n
f\left(\widehat{T}^\nu(x,0)\right)=\sum\limits_{\nu}^n\frac{(1-\alpha_\nu^\ast)\alpha_{\nu+1}}
{\alpha_{\nu+1}+\alpha_\nu^\ast},
$$
$$
\sum\limits_{\nu=1}^n
f\left(\widehat{T}^\nu(x,y)\right)=\sum\limits_{\nu}^n\frac{(1-\tilde{\alpha}_\nu)\alpha_{\nu+1}}
{\alpha_{\nu+1}+\tilde{\alpha}_\nu}.
$$

Рассмотрим ряд

$$
\sum\limits_{\nu=1}^\infty|f(\widehat{T}^\nu(x,y))-f(\widehat{T}^\nu(x,0))|=
\sum\limits_{\nu=1}^\infty\left|\frac{(1-\tilde{\alpha}_\nu)\alpha_{\nu+1}}
{\alpha_{\nu+1}+\tilde{\alpha_\nu}}
-\frac{(1-\alpha_\nu^\ast)\alpha_{\nu+1}}{\alpha_{\nu+1}+\alpha_\nu^\ast}\right|=
$$
$$
=\sum\limits_{\nu=1}^\infty\alpha_{\nu+1}\left|\frac{\tilde{\alpha}_\nu-\alpha_\nu^\ast+\tilde{\alpha_\nu}\alpha_{\nu+1}-
\alpha_\nu^\ast
\alpha_{\nu+1}}{(\alpha_{\nu+1}+\alpha_\nu^\ast)(\alpha_{\nu+1}+\tilde{\alpha}_\nu)}\right|\leqslant
\sum\limits_{\nu=1}^\infty\frac{1+\alpha_{\nu+1}}{\alpha_{\nu+1}}|\tilde{\alpha}_\nu-\alpha_\nu^\ast|\leqslant
$$
$$
\leqslant\sum\limits_{\nu=1}^\infty\frac{2y}{q_\nu(q_\nu+yp_\nu)}
\leqslant\sum\limits_{\nu=1}^\infty\frac{2}{q_\nu(q_\nu+p_\nu)}<4.
$$
Лемма 5 доказана.

{\bf Следствие 4.} {\it Для любого $y\in[0;1]$ выполняется
равенство}
$$
\lim\limits_{n\rightarrow\infty}\frac{1}{n}\sum\limits_{\nu=1}^n
f(\widehat{T}^\nu(x,y))=\lim\limits_{n\rightarrow\infty}\frac{1}{n}\sum\limits_{\nu=1}^n
f(\widehat{T}^\nu(x,0)).
$$

Обозначим через $R$ множество тех точек $(x,y)\in[0,1)^2$, для
которых выполняется равенство
$$
\lim\limits_{n\rightarrow\infty}\frac{1}{n}\sum\limits_{\nu=1}^n
f(\widehat{T}^\nu(x,y))=\frac12.
$$
Мера множества $R$ равна 1. Рассмотрим проекцию
$$
\tilde{R}=\{x\in[0,1)|\exists y: (x,y)\in R\}.
$$
Мера множества $\tilde{R}$ тоже будет равна 1. Согласно формуле
(\ref{DogrI}) и следствию 4 выполняется равенство
$$
\lim\limits_{t\rightarrow\infty}\frac{I_\alpha(t)}{N(\alpha,t)}=
\lim\limits_{n\rightarrow\infty}\frac{G_n(\alpha)}{n}=\frac{1}{n}\sum\limits_{\nu=1}^n
f(\widehat{T}^\nu(x,0))=\frac12.
$$
Пункт 1) теоремы 1 доказан.

Для доказательства пункта 2) теоремы 1 нам понадобится т. Леви (см
\cite{Khinchin}).

{\bf Теорема Леви.} {\it Для почти всех $\alpha\in(0;1)$ имеет
место следующее равенство}
$$
\lim\limits_{n\rightarrow\infty}\frac{\ln
q_n}{n}=\frac{\pi^2}{12\ln2}.
$$
Из теоремы Леви получаем что для почти всех $\alpha$ выполнено
равенство
$$
\lim\limits_{t \to \infty}\frac{\ln
t}{N(\alpha,t)}=\frac{\pi^2}{12\ln 12}.
$$
Пункт 2) теоремы 1 получается из пункта 1) теоремы 1 и последнего
равенства. Теорема 2 полностью доказана.

{\bf 7. Доказательство теоремы 4.} Снова рассматриваем функцию $f$
определенную в (\ref{fufu}). Пусть
$R_1\subset[0,1)^2\times[0,1)^2$ это множество для тех точек
$({z}_1, {z}_2) =((x_1,y_1),(x_2,y_2))$, для которых выполняется
следствие 3. Его мера равняется 1. Для $({z}_1, {z}_2) \in {R}_1$
рассмотрим только те значения $n$ для которых суммы
$$
\sum\limits_{\nu=1}^nf(\widehat{T}^\nu
z_1)-\sum\limits_{\nu=1}^nf(\widehat{T}^\nu z_2),
\,\,\,\,\,\,\,\,\,\,\,
 \sum\limits_{\nu=1}^{n+1}f(\widehat{T}^\nu
z_1)-\sum\limits_{\nu=1}^{n+1}f(\widehat{T}^\nu z_2)
$$
имеют различные знаки. Тогда, поскольку каждое слагаемое в суммах
не превосходит единицы, видим, что
$$
\left|\sum\limits_{\nu=1}^nf(\widehat{T}^\nu
z_1)-\sum\limits_{\nu=1}^nf(\widehat{T}^\nu z_2)\right| \le 2.
$$
Рассмотрим проекцию
$$\tilde{R}_1=\{(x_1,x_2)\in[0,1)^2 | \exists y_1 \exists y_2 :
(x_1,y_1,x_2,y_2)\in R_1\}.$$ Мера множества $\tilde{R}_1$ тоже
равна 1. Для $(\alpha, \beta) \in \tilde{R}_1$  и рассматриваемых
значений $n$ с учетом леммы 5 получаем
$$
\left| \int_1^{q_n}\psi_\alpha (t)dt - \int_1^{r_n}\psi_\beta
(t)dt \right|=\left| G_n(\alpha) - G_n(\beta)\right| =
\left|\sum\limits_{\nu=1}^nf\left(\widehat{T}^\nu
(\alpha,0)\right)-\sum\limits_{\nu=1}^nf\left(\widehat{T}^\nu
(\beta,0)\right)\right|\le 10.$$ Теорема 4 доказана.


\newpage


\begin{thebibliography}{99}
\bibitem{Kornfeld} И. П. Корнфельд, Я. Г. Синай, С. В. Фомин,
Эргодическая теория, М., Наука (1980)

\bibitem{Khinchin} А. Я. Хинчин, Цепные дроби, М., Физматлит,
(1960)

\bibitem{Ivanov} В. А. Иванов, О рациональных приближениях
действительных чисел, Математические заметки, {\bf т. 23, № 1}
(1978) стр. 3-26

\bibitem{nakada}
Hitoshi Nakada, Metrical Theory for a Class of Continued Fraction
Transformations and Their Natural Extensions, Tokyo J. Math., {\bf
4, no. 2}, (1981), стр. 399--426

\bibitem{nakada1}
H. Nakada, Sh. Ito, S. Tanaka, On the invariant measure for the
transformation associated with some real continued fraction, Keio
Engrg. Rep. {\bf 30, no. 13}, (1977), p. 159--175

\bibitem{MosKan}
I.D. Kan, N.G. Moshchevitin, Approximations to two real numbers,
Uniform Distribution Theory, {\bf 5, no. 2}, (2010), p. 79--86

\bibitem{Furs}
H. F\"{u}rstenberg, Y. Katznelson. D. Ornstein, The ergodic
theretical proof of Szemeredi's theorem, Bulletin (New series) of
the American Mathematical Society, $7:3$, (1982), p. 527--552

\bibitem{Hal}
G. Halasz, Remarks on the remainder in Birkhoff's ergodic theorem,
Acta Mathematica Academiae Scientiarum Hungarica {\bf 28 (3-4)}
(1976) p. 389--395

\bibitem{Cusick}
Thomas W. Cusick, Mary E. Flahive, The Markoff and Lagrange
spectra, Mathematica surveys and monographs {\bf 30}, (1943)

\bibitem{Kraaim}
Karma Dajani and Cor Kraaikamp, A Note on the Approximation by
Continued Fractions under an Extra Condition, New York Journal of
Mathematics, {\bf 3A}, (1998), p. 69--80

\end{thebibliography}
\end{document}